\documentclass[11pt,a4paper]{article}

\usepackage{epsf,epsfig,amsfonts,amsgen,amsmath,amstext,amsbsy,amsopn,amsthm
}
\usepackage{amsmath,times,mathptmx}
\usepackage{amsfonts,amsthm,amssymb}
\usepackage{amsfonts}
\usepackage{graphics}
\usepackage{latexsym,bm}
\usepackage{amsfonts,amsthm,amssymb,bbding}
\usepackage{indentfirst}
\usepackage{graphicx}
\usepackage{color}
\usepackage[colorlinks=true,anchorcolor=blue,filecolor=blue,linkcolor=blue,urlcolor=blue,citecolor=blue]{hyperref}
\usepackage{float}
\usepackage{tikz}
\usepackage{enumerate}
\newtheorem{thm}{Theorem}[section]

\newtheorem{lemma}{Lemma}
\newtheorem{false statement}{False statement}

\theoremstyle{definition}

\newtheorem{claim}{Claim}

\newtheorem{conj}{Conjecture}

\newtheorem{prob}{Problem}[section]

\begin{document}
	
	\title{Spectral conditions for forbidden subgraphs in bipartite graphs
		\footnote{Supported by
			 National Natural Science Foundation of China (No. 12271162) and Natural Science Foundation of Shanghai (No. 22ZR1416300).}}
	\author{{Yuan Ren, Jing Zhang, Zhiyuan Zhang}\thanks{Corresponding author. E-mail addresses: 18228578699@163.com (R. Yuan),  jsam0331@163.com (Z. Zhang).}\\
		{\footnotesize School of Mathematics, East China University of Science and Technology, Shanghai 200237, China}}
	\date{}
	
		\maketitle {\flushleft\large\bf Abstract:} 
	A graph $G$ is $H$-free, if it contains no $H$ as a subgraph. A graph is said to be \emph{$H$-minor free}, if it does not contain $H$ as a minor. In recent years, Nikiforov asked that what is the maximum spectral radius of an $H$-free graph of order $n$? In this paper, we consider about some Brualdi-Solheid-Tur\'{a}n type problems on bipartite graphs. In 2015, Zhai, Lin and Gong proved that if $G$ is a bipartite graph with order $n \geq 2k+2$ and $\rho(G)\geq \rho(K_{k,n-k})$, then $G$ contains a $C_{2k+2}$ unless $G \cong K_{k,n-k}$ [Linear Algebra Appl. 471 (2015)]. Firstly, we give a new and more simple proof for the above theorem. Secondly, we prove that if $G$ is a bipartite graph with order $n \geq 2k+2$ and $\rho(G)\geq \rho(K_{k,n-k})$, then $G$ contains all  $T_{2k+3}$ unless $G \cong K_{k,n-k}$. Finally, we prove that among all outerplanar bipartite graphs on $n>344569$ vertices, $K_{1,n-1}$ attains the maximum spectral radius.

	\vspace{0.1cm}
	\begin{flushleft}
	\end{flushleft}

	\section{Introduction}

	A graph $G$ is called \emph{$H$-free}, if it does not contain a subgraph isomorphic to $H$, where the graph $H$ can also been called \emph{the forbidden subgraph}. The classic Tur\'{a}n's problem asks what is the maximum edge number of an $H$-free graph of order $n$, where the answer is known as the Tur\'{a}n number of $H$ and denoted by $ex(n, H)$. The study of Tur\'{a}n's problem was initiated by Mantel \cite {M07}, who showed that $ex(n, K_3)\leq \lfloor n^2/4\rfloor$, where $K_j$ means a $j$-clique. Let $A(G)$ be the adjacency matrix of a graph $G$ and \emph{the spectral radius} $\rho(G)$ be the largest eigenvalue of $A(G)$. In 2010, Nikiforov \cite{N10} raised the following \emph{Brualdi-Solheid-Tur\'{a}n} type problem:
	
	\begin{prob}[\cite{N10}]\label{pr1}
		Given a graph $H$ or a family $\mathbb{H}$, what is the maximum spectral radius of $H$-free($\mathbb{H}$-free) graph of order $n$?
	\end{prob}
	
	Much attention has been paid to Brualdi-Solheid-Tur\'{a}n type problems for various types of $H$, such as the clique \cite{BN07,W86}, the complete bipartite graph \cite{BG09,NI4}, the path \cite{N10}, the friendship graph \cite{CFTZ20}, the odd wheel \cite{CDTar21} and consecutive cycles \cite{GH19,LNar,NP20,N08,ZLar,ZZar}. For more details, we refer the reader to the survey \cite{N11}. In 2010, Nikiforov \cite{N10} raised the following Conjecture \ref{conj1} based on the $\{C_{2k+1}, C_{2k+2}\}$-free graphs and $C_{2k+2}$-free graphs, and also raised Conjecture \ref{conj2} based on $T_{2k+2}$-free graphs and $T_{2k+3}$-free graphs in the same paper, where $C_j$ and $T_j$ means the cycle and some tree of order $j$, respectively. Denote by $\vee$ and $\cup$ the join and union products, respectively. Define $S_{n,k}=K_k\vee I_{n-k}$ and $S_{n,k}^+=K_{k}\vee(I_{n-k-2}\cup K_2)$, where $I_j$ means an independent set of order $j$.

	\begin{conj}[\cite{N10}]\label{conj1}
		Let $k\geq 2$ and $G$ be a graph of sufficiently large order $n$.
		
		\noindent(a) If $\rho(G)\geq\rho(S_{n,k})$, then $G$ contains a $C_{2k+1}$ or $C_{2k+2}$ unless $G\cong S_{n,k}$;
		
		\noindent(b) If $\rho(G)\geq\rho(S_{n,k}^+)$, then $G$ contains a $C_{2k+2}$ unless $G\cong S_{n,k}^+$. 
	\end{conj}
	
	\begin{conj}[\cite{N10}]\label{conj2}
		Let $k\geq 2$ and $G$ be a graph of sufficiently large order $n$. 
		
		\noindent(a) If $\rho(G)\geq \rho(S_{n,k})$, then $G$ contains all trees of order $2k+2$ unless $G\cong S_{n,k}$. 
		
		\noindent(b) If $\rho(G)\geq \rho(S_{n,k}^+)$, then $G$ contains all trees of order $2k+3$ unless $G\cong S_{n,k}^+$. 
	\end{conj} 
	Partial cases of Conjecture \ref{conj1} and \ref{conj2} were solved by Zhai and Wang \cite{ZW12}, Nikiforov \cite{N09}, Zhai and Lin \cite{ZL20} and other researchers such as \cite{HLWGL21,LBWar2109,LBWar2112}. Recently, Cioaba, Desai and Tait completely solved them \cite{CDTar2205,CDTar2206}. A graph is bipartite if and only if it has no odd cycles. Since the
	extremal graphs of the above two conjectures are both non-bipartite graph, it is natural to consider about 
	Brualdi-Solheid-Turán type problem for bipartite graphs. In \cite{ZLG15}, Zhai, Lin and Gong first give a spectral condition to guarantee the existence of even cycles in bipartite graphs. Denote by $K_{s,t}$ the complete bipartite graph with two partitions of order $s$ and $t$.
	
	\begin{thm}[\cite{ZLG15}]
		Let $k$ be a positive number, $n \geq 2k+2$ and $G$ be a bipartite graph of order $n$. 
		If $\rho(G)\geq \rho(K_{k, n-k})$, then $G$ contains a $C_{2k+2}$ unless $G \cong K_{k, n-k}$.\label{th1.1}
	\end{thm}
	
	In this paper, we shall give a new proof of Theorem \ref{th1.1} by local edge maximality and eigenvector entry. Also, we give a spectral condition to guarantee the existence of trees in bipartite graphs. 
	
	\begin{thm}\label{th4}
		Let $k\geq 2$ and $G$ be a bipartite graph of sufficiently large order $n$. If $\rho(G)\geq \rho(K_{k,n-k})$, then $G$ contains all trees of order $2k+3$ unless $G\cong K_{k,n-k}$. 
	\end{thm}
	
	Given two graphs $H$ and $G$, $H$ is a \emph{minor} of $G$ if $H$ can be obtained by means of a sequence of vertex deletions, edge deletions and edge contractions. A graph is said to be \emph{$H$-minor free}, if it does not contain $H$ as a minor. Also, we turn our focus to extremal spectral problem on $H$-minor free graph. Naturally, we consider the following problem.
	\begin{prob}
		Given a graph $H$ or a family $\mathbb{H}$, what is the maximum spectral radius of $H$-minor free($\mathbb{H}$-minor free) graph of order $n$?
	\end{prob}
	As we known that every planar graph is $\{K_{3,3}, K_5\}$-minor free and every outerplanar graph is $\{K_{2,3}, K_4\}$-minor free.
	The study of spectral extremal problems on planar and outerplanar graphs has a rich history. Denote by $P_{j}$ the path of order $j$. In 1990s, Cvetkovi\'{c} and Rowlinson \cite{CR90} conjectured that for any outerplanar graph $G$, $\rho(G)\leq \rho(K_1\vee P_{n-1})$  with equality if and only if $G\cong K_1\vee P_{n-1}$. Almost at the same time, Boots and Royle \cite{BR91} and Cao and Vince \cite{CV93} independently conjectured that for any planar graph $G$ of order $n\geq 9$, $\rho(G)\leq \rho(P_2\vee P_{n-2})$ with equality if and only if $G\cong P_2\vee P_{n-2}$. Subsequently, many scholars took an interest in these two conjectures (see \cite{CV93,H95,H98,SH00}). Ellingham and Zha \cite{EZ00} proved that $\rho(G)\leq 2+\sqrt{2n-6}$ for a planar graph $G$. Dvo\v{r}\'{a}k and Mohar \cite{DM10} proved that $\rho(G)\leq \sqrt{8\Delta-16}+3.47$ for a planar graph $G$ with the maximum degree $\Delta$. Based on these previous research, these two conjectures were confirmed by Tait and Tobin \cite{TT17} for sufficiently large $n$ in 2017, and especially, the conjecture on outerplanar graph was completely confirmed by Lin and Ning \cite{LN21} in 2021. Since the
	extremal graphs of the above two conjectures are also both non-bipartite graph, it is natural to consider about the problems for bipartite graphs.
	For planar bipartite graphs, Hong and Shu \cite{HS99} proved that $\lambda(G)\geq -\sqrt{2n-4}$ with equality if and only if $G\cong K_{2,n-2}$, where $\lambda(G)$ represents the least eigenvalue of  $A(G)$. Since that if $G$ is a bipartite graph, then we have $\rho(G)=-\lambda(G)$, thus we easily get the conclusion that among all planar bipartite graphs on $n$ vertices, $K_{2, n-2}$ attains the maximum spectral radius. For outerplanar bipartite graphs, we give the following extremal spectral result.
	
	\begin{thm}\label{th6}
		Among all outerplanar bipartite graphs on $n>344569$ vertices, $K_{1,n-1}$ attains the maximum spectral radius.
	\end{thm}
	
	\section{The proof of Theorem \ref{th1.1}}\label{sec3}
	
	In this section, we prove Theorem \ref{th1.1}. Before beginning our proof, we first give a lemma. We denote the edge set of $G$ by $E(G)$ and $e(G)=|E(G)|$.
	\begin{lemma}[\cite{ZLG15}]
	Let $G=\left\langle X, Y\right\rangle$ be a bipartite graph where $|X|\geq r$ and $|Y|\geq r-1\geq 1$. If $G$ does not contain a copy of $P_{2r+1}$ with both endpoints in $X$, then 
	$$e(G)\leq (r-1)|X|+r|Y|-r(r-1).$$
	Equality holds if and only if $G\cong K_{|X|, |Y|}$, where $|X|=r$ or $|Y|=r-1$. \label{lem1} 
	\end{lemma}
	In the following, we would like to introduce an eigenvector technique, \emph{Rayleigh quotient}. A Rayleigh quotient is a scalar of the form $y^\top Ay/y^\top y$ where $y$ is a non-zero vector in $\mathbb{R}^n$. The supremum of the set of such scalars is the largest
	eigenvalue $\rho$ of $A$, equivalently,
	\begin{equation*}
		\rho=\sup\{x^\top Ax: x\in\mathbb{R}^n, \Vert x\Vert=1\}.
	\end{equation*}
 	For a graph $G$ and any vertex $u\in V (G)$, we denote by $N_G(u)$ the neighborhood of $u$ in $G$. For subsets $X,Y\subset V(G)$, we write $E(X)$ for the set of edges induced by $X$ and $E(X,Y)$ for the set of edges with one endpoint in $X$ and another endpoint in $Y$. Define $e(X)=|E(X)|$ and $e(X,Y)=|E(X,Y)|$. Also, define $N_A(u) = N_{G^*}(u)\cap A$ and $d_A(u)=|N_A(u)|$ for any vertex $u\in V(G^*)$ and subset $A\subset V(G)$.
 	
	\noindent{\underline{\textbf{The proof of Theorem \ref{th1.1}.}}} Suppose that $G^*$ attains the maximum spectral radius among all $C_{2k+2}$-free bipartite graphs. We first assume that $G^*$ is disconnected. Let $G_1,...,G_\eta$ be all components of $G^*$ for $\eta\ge 2$, and let $G'$ be a connected graph obtained from $G^*$ by adding $\eta-1$ edges. By the Rayleigh quotient and the Perron-Frobenius theorem,  $\rho(G')> \rho(G^*)$. Also note that if $G^*$ is $C_{2k+2}$-free, then $G'$ is also $C_{2k+2}$-free. It is contrary to the maximality of $\rho(G^*)$. For this reason, in what follows, we always assume that $G^*$ is connected. Note that $K_{k,n-k}$ is $C_{2k+2}$-free. Then
	\begin{equation}
		\rho(G^*)\geq\rho(K_{k,n-k})=\sqrt{k(n-k)}.\label{2.1}
	\end{equation}
	Assume that $x=(x_1,x_2,\cdots,x_n)$ is the Perron vector of $G^*$ and $z$ is the vertex with the maximum entry of $x$. Let $A=N_{G^*}(z)$, $B=V(G^*)\setminus (A\cup \{z\})$ and
	\begin{equation}
		\gamma(z)=|A|+2e(A)+e(A,B).\label{2.2}
	\end{equation}
	Since $G^*$ is a bipartite graph, $e(A)=0$.  Note that
	$$\begin{aligned}
	\rho^2(G^*)x_z	& =  \rho(G^*)\sum\limits_{\substack{u\in N_{G^*}(z)}} x_u=\sum\limits_{\substack{u\in N_{G^*}(z)}}\sum\limits_{\substack{v\in N_{G^*}(u)}}x_v \\
		& =  |A|x_z +\sum\limits_{\substack{u\in A}}d_A(u)+\sum\limits_{\substack{v\in B}}d_A(v)x_v\\
		& \leq  |A|x_z +2e(A)x_z+e(A,B)x_z=\gamma(z)x_z.
	\end{aligned}$$

	Thus, combining this with (\ref{2.1}),
	\begin{equation}\label{2.3}
		k(n-k)\leq\rho^2(G^*)\leq \gamma(z).
	\end{equation}
	In the following, we prove a claim to finish our proof.
	\begin{claim}
		There exists a path $P_{2k+1}$ consisting of edges in $E(A,B)$ with both endpoints $u,v\in A$ containing in $G^*$ unless $G\cong K_{k,n-k}$.\label{claim2.1}
	\end{claim}

	\begin{proof}
		Otherwise, by Lemma \ref{lem1}, $$e(A,B)\leq (k-1)|A|+k|B|-k(k-1).$$ Then by (\ref{2.2}), $e(A)=0$ and $|A|+|B|+1=n$, we have
		$$\begin{aligned}
		\gamma(z)	& =  |A|+2e(A)+e(A,B) \\
			& \leq  k|A|+k|B|-k(k-1) \\
			& =  k(|A|+|B|+1-k)\\
			& =  k(n-k).
		\end{aligned}$$
		Combining this with (\ref{2.1}) and (\ref{2.3}), it is a contradiction unless the equality holds. Since the equality holds if and only if $G^*\cong K_{k,n-k}$ by Lemma \ref{lem1}, as desired.
	\end{proof}
	By Claim \ref{claim2.1}, since $u,v\in A$, $G^*$ contains a $C_{2k+2}:=zuP_{2k+1}v$ unless $G^*\cong K_{k,n-k}$, which is contrary to $C_{2k+2}$-free, we complete the proof.  $\hfill\square$
	
	In the proof of Theorem \ref{th1.1}, we aware that we can use the condition that $G$ is $C_3$-free instead of that $G$ is bipartite to finish the proof. We can obtain a general result as follows.
	\begin{thm}
		Let $k$ be a positive number, $n \geq 2k+2$ and $G$ be a $C_3$-free graph of order $n$. 
		If $\rho(G)\geq \rho(K_{k, n-k})$, then $G$ contains a $C_{2k+2}$ unless $G \cong K_{k, n-k}$.
	\end{thm}
	
		\section{The proof of Theorem \ref{th4}}\label{sec4}
	
	In this section, we will give the proof of Theorem \ref{th4}. Before beginning our proof, we first give some notations and lemmas. For $k \geq 2$, let $\mathbb{T}_k$ denote the set of all trees on $2k + 3$ vertices. We use $\mathbb{G}_{n,T}$ to denote the set of $T$-free bipartite graphs on $n$ vertices for $T\in \mathbb{T}_k$. Let $G_{n,T}$ be a graph with maximum spectral radius among all bipartite graphs in $\mathbb{G}_{n,T}$. 
	
	\begin{lemma}[\cite{CDTar2206}]\label{lemma 3.1}
		Let $y$ be a non-negative non-zero vector and $c$ be a positive constant. For every non-negative symmetric matrix $A$, if $Ay\geq cy$ entrywise, then $\rho(A)\geq c$.
	\end{lemma}
	
	\begin{proof}
		We have $y^{\top}Ay\geq y^{\top}cy$ since $Ay\geq cy$ entrywise. Then by Rayleigh quotient, $\rho(A)\geq \frac{y^{\top}Ay}{y^{\top}y}\geq c$ follows.
	\end{proof}
	
	The following lemma is an exercise in \cite{BM76}. For short, we use $H\subseteq G(H\subset G)$ to express that $H$ is a (proper) subgraph of $G$ below.
	
	\begin{lemma}[\cite{BM76}]\label{0}
		Let $T$ be an arbitrary tree on $k + 1$ vertices. If $G$ is a graph with its minimum degree $\delta(G)\geq k$, then $T\subseteq G$.
	\end{lemma}
	
		\begin{lemma}\label{6}
		If $G$ is a graph with size $e(G)> kn$, then $H\subseteq G$ and its minimun degree $\delta(H)\geq k+1$.
	\end{lemma}
	
	\begin{proof}
		Let $H$ be an induced subgraph of $G$  defined by  a sequence of graphs $G_1, G_2, \cdots, G_{t+1}$ such that:
		\begin{enumerate}[(1)]
			\item $G_1=G$, $G_{t+1}=H$;
			\item for every integer $j\in [1,t]$ (if $t\geq 1$), there exists one vertex $v_{i_j}\in V(G_j)$ such that $d_{G_j}(v_{i_j})< k+1$ and $G_{j+1}=G_j-\{v_{i_j}\}$;
			\item for every vertex $v\in V(G_{t+1})$ (if any), $d_{G_{t+1}}(v)\ge k+1$.
		\end{enumerate}
		By the construction, $|V(G)\setminus V(H)|=t$, and if $t\ge 1$ then $d_{G_j}(v_{i_{j}})\le k$ for $1\leq j\leq t$. We assert that $H$ is not a null graph. Otherwise, assume that $H$ is a null graph, i.e., $t=n$. Then $$e(G)=d_G(v_{i_1})+d_{G_2}(v_{i_2})+\cdots+d_{G_t}(v_{i_t})\le tk=kn,$$
		contrary to $e(G)>kn$. Similarly, we can prove that $|V(H)|\geq k+2$ since $$e(G)=d_G(v_{i_1})+d_{G_2}(v_{i_2})+\cdots+d_{G_t}(v_{i_t})+e(H)\le tk+(n-t)(n-t-1)\leq tk+(n-t)k= kn,$$ which is contrary to $e(G)>kn$ when $|V(H)|=n-t\leq k+1$. By the construction of $H$, $\delta(H)\geq k+1$, as desired.
	\end{proof}

	\begin{lemma}\label{lemma 3.2+3.3}
		Let $n$ be a sufficiently large integer, for any tree $T_t$, we have
		
		(\textsl{a}) $ ex(n, T_t)\leq (t-2)n.$
		
		(\textsl{b})  $T_t\subset K_{\lfloor\frac{t}{2}\rfloor,t-1}$. 
	\end{lemma}
	
	\begin{proof}
		We first give the proof of (\textsl{a}). Suppose to the contrary that $ ex(n, T_t)> (t-2)n.$ Let $G$ be a graph with size $e(G)=ex(n,T_t)>(t-2)n$, By Lemma \ref{6}, if $e(G)>(t-2)n$, then there exists a subgraph $H\subseteq G$ with minimum degree $\delta(H)\geq t-1$. Thus, $H$ contains $T_t$ by Lemma \ref{0}, a contradiction. This implies $ex(n, T_t)\leq (t-2)n$.
		
		Then, we give the proof of (\textsl{b}). Note that every tree is bipartite. Then we can divide the vertex set of $T_t$ into two partitions. Since the larger partition of $T_t$ has $a\leq t-1$ vertices and the smaller partition has $b\leq \lfloor\frac{t}{2}\rfloor$ vertices, $T_t\subseteq K_{a,b}\subseteq K_{t-1,\lfloor\frac{t}{2}\rfloor}$, the proof is complete.
	\end{proof}

	For a graph $G$ and $v\in V(G)$, let $N_i(v)$ be the vertices at distance $i$ from $v$ and $d_i(v)=|N_i(v)|$. Specially, $d_1(v)=d_G(v)$. For short, we use $\rho$ and $d(v)$ instead of $\rho(G)$ and $d_G(v)$ if $G$ is clear from the context respectively.
	 
	\begin{lemma}\label{lemma 3.4}
	For two integers $k\geq 2$ and $n\geq k+2$,
	\begin{equation*}\sqrt{k(n-k)}\leq \rho(G_{n,T})\leq \sqrt{(2k+1)n}\leq \sqrt{3kn}.
	\end{equation*}
	\end{lemma}

	\begin{proof}
	Note that $K_{k,n-k}\in \mathbb{G}_{n,T}$ and $G_{n,T}$ is the graph with the maximum spectral radius in $\mathbb{G}_{n,T}$. Then $\rho(G_{n,T})>\rho(K_{k,n-k})=\sqrt{k(n-k)}$. Since $G_{n,T}$ is a bipartite graph, $e(N_1(v))=0$. By eigenvalue-eigenvector equation, we have
	\begin{center}\begin{equation}\label{8}
			\begin{aligned}
				\rho^2 & \leq \max_{v\in V} \left\{{\sum_{u\in N_{G_{n,T}}(v)}{A^2_{v,u}}}\right\} =\max_{v\in V}\left\{{\sum_{u\in N_{G_{n,T}}(v)}{d(u)}}\right\}\\
				&=\max_{v\in V}\{{d(v)+2e(N_1(v))+e(N_1(v),N_2(v))}\}\\
				&=\max_{v\in V}\{{d(v)+e(N_1(v),N_2(v))}\}.\end{aligned}
\end{equation}\end{center}
	By Lemma \ref{lemma 3.2+3.3} (\textsl{a}), $e(G_{n,T})\leq(2k+1)n$. Combining this with (\ref{8}), we have
	\begin{center}\begin{equation*}
			\begin{aligned}
				\rho^2 \leq e(G_{n,T})\leq (2k+1)n\leq 3kn,\end{aligned}
	\end{equation*}\end{center}
	as desired.
	\end{proof}
	Let $\eta$ be a positive constant which only depends on $k$:

	\begin{equation}\label{i.e2}
	\eta<\min\left\{\frac{1}{10k}, \frac{1}{(2k+2)(16k^2)}\cdot\big(16k^2-1-\frac{4(16k^3-1)}{5k}\big), \frac{32k^2-1}{64k^3+8k}\right\},
	\end{equation}
	Let $x=(x_1,\cdots,x_n)$ be the eigenvector corresponding to the eigenvalue $\rho(G_{n,T})$  such that $x_z=\max_{v\in V(G)} x_v =1.$ Then we introduce two positive constants that depend on $\eta$ and $k$:
	\begin{equation}\label{i.e3}
	\epsilon<\min\left\{\frac{\eta}{2}, \frac{1}{16k^3}, \frac{1}{8k}, \frac{\eta}{32k^3+2}\right\}\text{, }
\end{equation}
\begin{equation}\label{i.e5}
	\alpha<\min\left\{\eta, \frac{\epsilon^2}{3k}\right\}.
	\end{equation}
	For convenience, we individually introduce the following inequality, which can be easily deduced by (\ref{i.e2}) and (\ref{i.e3}):
\begin{equation}\label{ie2.1.5}
	\eta-2\varepsilon>(32k^3+2)\varepsilon-2\cdot\varepsilon={32k^3}\varepsilon.
\end{equation}
	In the following, we define two vertex subsets of $G_{n,T}$ that depend on $\alpha$ and $x$:
	$$L:=\{v\in V|x_v\geq \alpha\},~~M:=\{v\in V|x_v\geq \alpha/3\}.$$ 
	
	\begin{lemma}\label{lemma 3.5}
		For $n\geq\frac{25k^3}{9k^2-16k-4}$, we have
		$|L|\leq \frac{5\sqrt{kn}}{\alpha}$ and $|M|\leq  \frac{15\sqrt{kn}}{\alpha}$.
	\end{lemma}
	
	\begin{proof}
		For any $v\in V(G_{n,T})$, we have that
		\begin{equation}\label{ie3.5.4}
		\rho x_v=\sum_{u\sim v}{x_u}\leq d(v).
		\end{equation}
		Combining this with the lower bound in Lemma \ref{lemma 3.4}, we obtain $\sqrt{k(n-k)}x_v\leq \rho x_v\leq d(v)$. Summing over all vertices $v\in L$ gives
		\begin{equation*}
		|L|\sqrt{k(n-k)}\alpha \leq \sum_{v\in L}{d(v)}\leq \sum_{v\in V}{d(v)}\leq 2e(G_{n,T})\leq 2\max_{T\in \mathbb{T}_k}{ex(n,T)}\leq (4k+2)n,
		\end{equation*}
		where the last inequality holds by Lemma \ref{lemma 3.2+3.3} (a). Thus, we obtain that $|L|\leq \frac{5\sqrt{kn}}{\alpha}$ for $n\geq\frac{25k^3}{9k^2-16k-4}$. By a similar analysis, we have 
		$$|M|\sqrt{k(n-k)}\cdot\frac{\alpha}{3} \leq (4k+2)n,$$ which follows that $|M|\leq \frac{15\sqrt{kn}}{\alpha}$.
	\end{proof}
		Now we give the proof of Theorem \ref{th4}.
		
	\noindent{\underline{\textbf{The proof of Theorem \ref{th4}.}}} First of all, we assert that the graph $G_{n,T}$ is connected.
	If not, $G_{n,T}$ has $r\geq 2$ components: $C_1, \dots ,C_r$. Let $C_1$ be a component with $\rho(G_{n,T})=\rho(C_1)$. Note that $x$ is an eigenvector responding to $\rho(G_{n,T})$ and $x_z=1$. Without loss of generality, we assume that $z\in C_1$. By (\ref{ie3.5.4}),  
	$$d(z)\geq \rho(G_{n,T})\geq \sqrt{k(n-k)}.$$ 
	Since other components have at least one vertex, there exists a vertex $u\in C_s, s\ne 1$. Let 
	$$\hat G_{n,T}=G_{n,T}-\{uv:\text{for all } v\in N_1(u)\}+\{uz\}.$$ 
	Let $C_1^{'}=C_1+\{uz\}$. It's clear that $C_1\subset C_1^{'}\subseteq \hat{G}_{n,T}$. Then $$\rho(\hat G_{n,T})\geq \rho(C_1^{'})>\rho(C_1)=\rho(G_{n,T}).$$ Moreover, we claim that $\hat G_{n,T}$ contains no tree in $\mathbb{T}_k$. If not, we suppose $\hat{G}_{n,T}$ contains $T$. By the construction of $\hat G_{n,T}$, we obtain that $u$ must be a leaf of $T$ and is adjacent only to the vertex $z$ in $\hat{G}_{n,k}$. Note that $d(z)\geq \sqrt{k(n-k)}$. Thus, the vertex $z$ is adjacent to at least $\sqrt{k(n-k)}-(2k+1)$ vertices in $V(C_1)\setminus V(T)$, and so there exists $w\in N_1(z)\setminus V(T)$, such that $$T'=T-uz+zw\cong T\text{, }T'\subseteq C_1\subset G_{n,T}.$$ It is impossible because $G_{n,T}$ contains no $T\in\mathbb{T}_k$. Since $\rho(\hat G_{n,T})>\rho(G_{n,T})$ and $\hat G_{n,T}$ contains no tree in $\mathbb{T}_k$, it is contrary to the maximality of $G_{n,T}.$ 
	
%
%

	Let $L_i(v)=L\cap N_i(v)$ and $M_i(v)=M\cap N_i(v)$. For short, we will use $N_i$, $L_i$ and $M_i$ instead of $N_i(v)$, $L_i(v)$ and $M_i(v)$.
	\begin{claim}\label{claim4.1}
		For any vertex $v\in L$ with $x_v=c\geq\alpha$, we have $d(v)\geq\frac{n\alpha}{8k}$. 
	\end{claim}
	
	\begin{proof}
		Suppose to the contrary that there exists a vertex $v\in L$ such that
		\begin{equation}\label{con1}
		d(v)<\frac{n\alpha}{8k}.
		\end{equation}
		Combining eigenvalue-eigenvector equation on some vertex $v$ with Lemma \ref{lemma 3.4} and $e(N_1)=0$ as $G_{n,T}$ is a bipartite graph, we have 
		\begin{equation}\label{ie4.1.1}
		\begin{aligned}
		k(n-k)c&\leq \rho^{2}c\\&=\rho^{2}x_{v}\\
		&=\sum\limits_{\substack{u\in N_1(v)}}\sum\limits_{\substack{w\in N_1(u)}}x_{w} \\
		 &\leq d(v)c+2e(N_1(v))+\sum\limits_{\substack{u\in N_1(v) }}\sum\limits_{\substack{w\in N_1(u),\\w\in N_{2}(v)}}x_{w}\\
		&=d(v)c+\sum\limits_{\substack{u\in N_1(v)}}\sum\limits_{\substack{w\in N_1(u),\\w\in N_{2}(v)}}x_{w}\\
		&\leq d(v)+\sum\limits_{\substack{u\in N_1(v)}}\sum\limits_{\substack{w\in N_1(u),\\w\in M_{2}(v)}}x_{w}+\sum\limits_{\substack{u\in N_1(v)}}\sum\limits_{\substack{w\sim u,\\w\in N_{2} \setminus M_{2}(v)}}x_{w} .
		\end{aligned}
		\end{equation}
		Combining this with Lemma \ref{lemma 3.5} and (\ref{con1}), we have
		$$
		\begin{aligned}
		\sum\limits_{\substack{u\sim v}}\sum\limits_{\substack{w\sim u,\\w\in M_{2}(v)}}x_{w}&\leq e(N_{1}(v),M_{2}(v))\\ &\leq ex(d(v)+|M|,T)& (\text{since }G_{n,T}[N_1, M_2]\text{ is a }T\text{-free bipartite graph.})\\
		&\leq (2k+1)(d(v)+|M|)&(\text{by Lemma \ref{lemma 3.2+3.3}})\\&<(2k+1)\left(\frac{n\alpha}{8k}+\frac{15\sqrt{kn}}{\alpha}\right).
		\end{aligned}$$
		Furthermore, for $n\geq\left(\frac{15\cdot24(2k+1)k\sqrt{k} }{(2k-3)\alpha^2}\right)^2$, we have
		
		$$(2k+1)\left(\frac{n\alpha}{8k}+\frac{15\sqrt{kn}}{\alpha}\right)\leq (2k+1)n\left(\frac{\alpha}{8k}+\frac{\alpha(2k-3)}{24k(2k+1)}\right)=\frac{n\alpha}{3}.$$
		Combining this with Lemma \ref{lemma 3.2+3.3} (\textsl{a}) and (\ref{ie4.1.1}), we get
		$$k(n-k)c\leq\frac{n\alpha}{8k}+\frac{n\alpha}{3}+e(N_{1}(v),N_{2} \backslash M_{2}(v))\frac{\alpha}{3}
		\leq \frac{n\alpha}{8k}+\frac{n\alpha}{3}+(2k+1)n\frac{\alpha}{3} = (\frac{1}{8k}+\frac{2k}{3}+\frac{2}{3})n\alpha,$$
		which implies a contradiction for $n>\frac{24k^3}{8k^2-16k-3}$ and $x_v=c\geq\alpha$ since $v\in L$, as desired.
		
	\end{proof}
	Again by Lemma \ref{lemma 3.2+3.3} (\textsl{a}), we have
	$$(2k+1)n\geq e(G_{n,T})\geq \frac{1}{2}\sum\limits_{\substack{v\in L}}d(v)\geq |L|\frac{n\alpha}{16k}.$$
	Then we get the refinement of the size of set $L$ that 
	\begin{equation}\label{size of L}
	|L|\leq \frac{16k(2k+1)}{\alpha}.
	\end{equation}
	Next we will refine the degree estimates of the vertices in $L^{'}:=\{v\in L|x_v\geq \eta\}$. Let $S:=V\setminus L=\{v\in V|x_v<\alpha\}$ and $S_i(v):=S\cap N_i(v)$. For short, $S_i=S_i(v)$ if $v$ is clear from the context.
	
	\begin{claim}\label{claim4.2}
	Let $v$ be a vertex in $L^{'}$ with Perron weight $x_{v}=c$. Then $d(v)\geq cn-\epsilon n.$
	\end{claim}
	
	\begin{proof}
	 Suppose to the contrary that 
	 \begin{equation}\label{con2}
	 d(v)<cn-\epsilon n. 
	 \end{equation}
	 Combining eigenvalue-eigenvector equation on the vertex $v$ with $e(S_1)=e(L_1,S_1)=0$ as $G_{n,T}$ is bipartite,  we have 
		
		\begin{equation}\label{ie4.1.3}
		\begin{aligned}
		k(n-k)c&\leq\rho^{2}c=\sum\limits_{\substack{u\in N_1( v)}}\sum\limits_{\substack{w\in N_1(u)}}x_{w}  = d(v)c+\sum\limits_{\substack{u\in N_1(v)}}\sum\limits_{\substack{w\in N_1(u),\\w\neq v}}x_{w}\\
		&\leq d(v)c+\sum\limits_{\substack{u\in S_{1}}}\sum\limits_{\substack{w\in N_1(u),\\w\in L_{2}}}x_{w}+2e(S_1)\alpha+2e(L)+e(L_1,S_1)\alpha+e(N_{1},S_{2})\alpha .\\
		&\leq d(v)c+\sum\limits_{\substack{u\in S_{1}}}\sum\limits_{\substack{w\in N_1(u),\\w\in L_{2}}}x_{w}+2e(L)+e(N_{1},S_{2})\alpha.
		\end{aligned}
		\end{equation}
		By Lemma \ref{lemma 3.2+3.3}, $e(N_{1},S_{2})\leq (2k+1)n.$
		Combining this with (\ref{size of L}), for $n\geq \frac{16^2k^2(2k+1)^2}{(k-1)\alpha^3}$ we deduce that\\
		$$2e(L)+e(N_{1},S_{2})\alpha \leq 2{|L|\choose 2}+(2k+1)n\alpha \leq |L|^2+(2k+1)n\alpha\leq3kn\alpha. $$
		Thus, by (\ref{i.e5}) and (\ref{ie4.1.3}), we obtain that\\  	
		\begin{equation}\label{ie4.1.4}
		\begin{aligned}
		k(n-k)c&\leq\rho^2c\\&\leq d(v)c+\sum\limits_{\substack{u\in S_{1}}}\sum\limits_{\substack{w\sim u,\\w\in L_{2}}}x_{w}+3kn\alpha \\&\leq d(v)c+e(S_{1},L_{2})+3kn\alpha \\&\leq d(v)c+e(S_{1},L_{2})+\epsilon^{2}n.
		\end{aligned}
		\end{equation}
	Note that $v\in L^{'}$. Then $\eta\geq x_v=c$. Thus, by (\ref{i.e3}), we have $$\varepsilon<\eta\leq c\leq 1.$$ Combining this with (\ref{con2}) and (\ref{ie4.1.4}), we deduce that
				$$(k-c+\epsilon)nc-k^2c=k(n-k)c-(cn-\epsilon n)c\leq (k(n-k)-d(v))c\leq e(S_{1},L_{2})+\epsilon^{2}n,$$
				which follows that
		\begin{equation}\label{ie4.1.5}
		e(S_{1},L_{2})\geq (k-c)nc+\epsilon nc-\epsilon ^{2}n-k^2c\geq (k-1)nc+\epsilon^2n-k^2c.
		\end{equation}
				
		Next we will prove that there are at least $\delta n$ vertices in $S_1$ with degree at least $k$ in  $G_{n,T}[S_1, L_2]$, where $\delta:=\frac{\epsilon\alpha}{16k(2k+1)}$. If not, there are at most $\delta n$ vertices in $S_1$ with degree at least $k$ in $G_{n,T}[S_1, L_2]$. Thus, by $|S_1| \leq d(v)< cn-\epsilon n$ and (\ref{size of L}), $e(S_1, L_2) < (k-1)|S_1| + |L|\delta n \leq (k-1)(c-\epsilon)n + \epsilon n$. Combining this with (\ref{ie4.1.5}) gives $-(k-2)n\epsilon> \epsilon^2n-k^2c\geq0$ for $n\geq \frac{k^2}{\varepsilon^2}$, which is impossible because $k\geq 2$ and $\varepsilon>0$.
		
		Let $B$ be the subset of $S_1$ in which each vertex has degree at least $k$ in $G_{n,T}[S_1, L_2]$. Thus $|B| \geq \delta n.$ Since there are only ${|L_2|\choose k}\leq{|L|\choose k}\leq{16k(2k+1)/\alpha\choose k} $ options for every vertex in $B$ to choose a set of $k$ neighbors from, there exists some set of $k$ vertices in $L_2$ with at least $\delta n/{|L|\choose k}\geq \delta n/{16k(2k+1)/\alpha\choose k} =\frac{\epsilon \alpha}{16k(2k+1)}n/{16k(2k+1)/\alpha\choose k}\geq 2k+2$ common neighbors in $B$, and so $K_{k+1,2k+2}\subset G_{n,T}[S_1, L_2]$, which is contrary to Lemma \ref{lemma 3.2+3.3}(b) that $G_{n,T}$ contains $T$ as desired.
	\end{proof} 
	
	\begin{claim}\label{claim4.3}
	 $(1-\epsilon)kn\leq e(S_{1},\{z\}\cup L_{2})\leq (k+\epsilon)n$.
	\end{claim}
	
	\begin{proof}
		Recall that $x_z=1$. Since $e(S_{1},\left\lbrace z\right\rbrace \cup L_{2})=e(S_{1},L_{2})+d(z)-|L_1|$, by (\ref{ie4.1.4}), we have
		\begin{equation}\label{i.e4}
		\begin{aligned}
		k(n-k)\cdot 1&\leq d(z)+e(S_{1}, L_{2})+\epsilon ^{2}n\\&=e(S_{1},\{z\}\cup L_{2})+|L_1|+\epsilon^2n
		\\&\leq e(S_{1},\{z\}\cup L_{2})+|L|+\epsilon^2n.
		\end{aligned}
		\end{equation}
		By (\ref{i.e3}), $\varepsilon<\frac{1}{26k^3}$. Combining this with (\ref{size of L}), $e(S_1,\{z\cup L_2\})\geq (1-\varepsilon)kn$. For the upper bound of Claim \ref{claim4.4}, suppose to the contrary that $e(S_{1},\left\lbrace z\right\rbrace \cup L_{2})>(k+\epsilon)n$. 
		We first assert that there are at least $\delta n$ vertices in $S_{1}$ with degree at least $k$ in $G_{n,T}[S_{1},L_{2}]$, where $\delta:=\frac{\epsilon\alpha}{16k(2k+1)}$. Otherwise by $|S_{1}|\leq n$ and (\ref{size of L}), we obtain that $e(S_{1},L_{2})<(k-1)|S_{1}|+|L|\delta n\leq (k-1)n+\epsilon n=(k+\epsilon-1)n$, which leads to a contradiction that $e(S_{1},\left\lbrace z\right\rbrace \cup L_{2})>(k+\epsilon)n$.  Therefore, let $D$ be the subset of $S_{1}$ with at least $\delta n$ vertices such that every vertex in $D$ has degree at least $k$ in $G_{n,T}[S_{1},L_{2}]$. For every vertex in $D$, there are at most ${|L|\choose k}\leq{16k(2k+1)/\alpha\choose k}$ options to choose a set of $k$ neighbors from. Thus, we have that there exists some set of $k$ vertices in $L_{2}\backslash \left\lbrace z\right\rbrace $ having a common neighborhood with at least $\delta n/{|L|\choose k}\geq \delta n/{16k(2k+1)/\alpha\choose k}  =\frac{\epsilon \alpha} {16k(2k+1)}n/{16k(2k+1)/\alpha\choose k}  \geq 2k+2$ vertices for $n$ large enough. So, $K_{k+1,2k+2}\subset G_{n,T}[S_{1},L_{2}\cup \left\lbrace z\right\rbrace ]$, which is contrary to Lemma \ref{lemma 3.2+3.3} (\textsl{b}), as desired.
	\end{proof}
	
	\begin{claim}\label{claim4.4}
		Let $v$ be a vertex with $v\in L^{'}:=\{v\in V|x_v\geq \eta\}$. We have $x_{v}\geq (1-\frac{1}{16k^{3}})$ and $d(v)\geq (1-\frac{1}{8k^{3}})n$.
	\end{claim}
	
	\begin{proof}
		First we prove that $x_{v}\geq (1-\frac{1}{16k^{3}})$ for any $v\in L^{'}$ through a contradiction by assuming that there is some vertex in $v\in L^{'}$ and $x_{v}< (1-\frac{1}{16k^{3}})$. Then refining (\ref{ie4.1.4}) with respect to the vertex $z$ we have that
		\begin{equation*}
		\begin{aligned}
		k(n-k)\cdot 1&\leq\rho^{2}\cdot 1\\&<e(S_{1}(z),L_{2}(z)\backslash \left\lbrace v\right\rbrace )+|N_{1}(z)\cap N_{1}(v)|x_{v}+\epsilon^{2}n\\
		&=e(S_1(z),L_2(z)\cup\{z\})-e(S_1(z),\{z\})-e(S_1(z),\{v\})+|N_{1}(z)\cap N_{1}(v)|x_{v}+\epsilon^{2}n\\
		&<(k+\epsilon)n-|S_{1}(z)\cap N_{1}(v)|+|N_{1}(z)\cap N_{1}(v)|\left( 1-\frac{1}{16k^{3}}\right)+\epsilon^{2}n\\
		&=kn+\epsilon n+|L_{1}(z)\cap N_{1}(v)|-|N_{1}(z)\cap N_{1}(v)|\frac{1}{16k^{3}}+\epsilon^{2}n.
		\end{aligned}
		\end{equation*}
		Thus, for sufficiently large $n$, we have $$|N_{1}(z)\cap N_{1}(v)|\frac{1}{16k^{3}}<\epsilon n+|L_1(z)\cap N_1(v)|+\epsilon^2n+k^2<\epsilon n+\epsilon^{2}n+|L|+k^2\leq 2\epsilon n,$$ by (\ref{size of L}). Recall that $|N_1(z)|\geq(1-\epsilon)n$ by Claim \ref{claim4.2}. Then there are at most $\epsilon n$ vertices that not contained in $N_1(z)$. Since $v\in L^{'}$, we have $x_{v}\geq \eta$  and  $d(v)\geq (\eta -\epsilon)n$ from Claim \ref{claim4.2}, and so $|N_{1}(z)\cap N_{1}(v)|\geq(\eta-\epsilon)n-\epsilon n= (\eta -2\epsilon)n>32k^{3}\epsilon n$ by (\ref{ie2.1.5}), a contradiction. This implies that $x_u>1-\frac{1}{16k^3}$. Combining this with Claim \ref{claim4.2} and (\ref{i.e3}), we have that $$d(v)\geq(x_v-\epsilon) n\geq(1-\frac{1}{16k^3}-\epsilon)n\geq (1-\frac{1}{8k^{3}})n,$$
		as desired.
	\end{proof}
	
	Now  we shall show $|L^{'}|= k$ by proving $|L^{'}|\geq k+1$ and $|L^{'}|\leq k-1$ is impossible. If $L^{'}\geq k+1$, then $G_{n,T}[S_{1},L^{'}]$ contains a $K_{k+1,2k+2}$ since the $(k+1)$-vertex set $B\subset L^{'}$ has a common neighbourhood with at least $(1-\frac{k+1}{8k^3})n>2k+2$ vertices for $n$ large enough, contrary to Lemma \ref{lemma 3.2+3.3}(\textsl{b}).
	Next, if $|L^{'}|\leq k-1$, then by (\ref{i.e4}), we have
	$$k(n-k)\leq \rho^{2}\leq e(S_{1},L^{'})+e(S_{1},L_{2} \backslash L^{'})\eta +\epsilon^{2}n\leq (k-1)n+(k+\epsilon)n\eta +\epsilon^{2}n<k(n-k),$$
	where the last inequality holds for $n$ large enough. This gives the contradiction. Thus 
	\begin{equation*}
	|L^{'}|=k.
	\end{equation*} 
	Now we have the above equation and every vertex in $L^{'}$ has degree at least $(1-\frac{1}{8k^3})n$. Thus, the common neighborhood $R$ of vertices in $L^{'}$ has at least $(1-\frac{1}{8k^2})n$ vertices. Let $E=V(G)\setminus(R\cup L^{'})$. Clearly, $|E|\leq \frac{n}{8k^2}$.

	\begin{claim}\label{claim4.5}
		For any vertex $v\in V(G_{n,T})$, we have $\rho x_v=\sum_{w\sim v}{x_w}\geq k-\frac{1}{16k^2}$.
	\end{claim}
	
	\begin{proof}
		We will divide the examinations of the Perron weight in the neighborhood of $v$ into three cases.
		
		\noindent\underline{\textbf{Case (a): }$v\in L^{'}$.}~ By Claim \ref{claim4.4}, $$\sum_{w\sim v}{x_w}=\rho x_v\geq \rho(1-\frac{1}{16k^3}).$$Combining this with $\rho\geq\sqrt{k(n-k)}>k$, $$\sum_{w\sim v}{x_w}\geq \rho(1-\frac{1}{16k^3})\geq k-\frac{1}{16k^2},$$
	 since $n$ is large enough.
		
		\noindent\underline{\textbf{Case (b): }$v\in R$.}~ Since $|L^{'}|=k$ and the vertex $v$ is a common neighbor of $L^{'}$, we deduce that $$\sum_{w\sim v}{x_w}\geq \sum_{\substack{w\sim v\\w\in L^{'}}}{x_w}\geq k(1-\frac{1}{16k^3})=k-\frac{1}{16k^2}.$$

		\noindent\underline{\textbf{Case (c): }$v\in E$.}~ Suppose to the contrary that $\sum_{w\sim u}{x_w}<k-\frac{1}{16k^2}$, given the graph $$H=G_{n,T}-\{vw:\text{for all } w\in N_1(v)\}+\{vu: \text{for all } u\in L^{'}\}.$$ Now since $\sum_{w\sim v}{x_w}<k-\frac{1}{16k^2}$, we have that $x^{\top}A(H)x>x^{\top}A(G_{n,T})x=\rho(G_{n,T})$ since $x$ is an eigenvector with respect to $\rho(G_{n,T})$, and $\rho(H)>\rho(G_{n,T})$ by Rayleigh quotient. However, we assert $H$ is $T$-free. Otherwise, if $H$ contains $T$, then $T$ contains the vertex $v$ since $G_{n,T}$ is $T$-free. Clearly, there exists some edge $e\in E(H)\cap E(T)$ which must be incident with the vertex $v$ and the other vertex $u\in L^{'}$. Now since $R$ has at least $(1-\frac{1}{8k^2})n>2k+2$ vertices, there must exist another vertex $v'\in R, v'\notin V(T)$ satisfies that for all $vu\in E(H)\cap E(T)$, $v'u\in E(G_{n,T})$. Thus, $G_{n,T}$ contains $T$, a contradiction, and so $H$ is $T$-free. It is contrary to the maximality with $G_{n,T}$, as desired.
		
		To sum up, we finish our proof of Claim \ref{claim4.5}.
		
	\end{proof}
	
	Since the graph $G_{n,T}$ is bipartite, the set $L^{'}$ must be contained in one partition, while $R$ is contained in the other set since $R$ is the common neighborhood of $L^{'}$. It follows that $e(R)=0$. Moreover, any vertex in $E$ is adjacent to at most $2k+1$ vertices in $R$, else $K_{k+1,2k+2}\subset G_{n,T}$ , which is contrary to Lemma \ref{lemma 3.2+3.3} (\textsl{b}). Finally, any vertex in $E$ is adjacent to at most $k-1$ vertices in $L^{'}$ by the definition of $E$.  Next we  prove that the vertex set $E=\emptyset$ for $n>\frac{128k}{3}$. 
	
	Assume to the contrary that $E\neq \emptyset$. Recall that any vertex $r\in R$ satisfies $x_r<\eta $. Therefore, any vertex $v\in E$ must satisfy
	$$
	\rho x_v=\sum\limits_{u\sim v}{x_u}=\sum\limits_{\substack{u\sim v\\u\in L^{'}\cup R}}{x_u}+\sum\limits_{\substack{u\sim v\\u\in E}}{x_u}< k-1+(2k+2)\eta +\sum\limits_{\substack{u\sim v\\u\in E}}{x_u}.
	$$
	Combining this with Claim \ref{claim4.5} and (\ref{i.e2}), we have
	\begin{equation*}
	\frac{\sum_{\substack{u\sim v\\u\in E}}{x_u}}{\rho x_v}> \frac{\rho x_v-(k-1)-(2k+2)\eta}{\rho x_v}\geq 1-\frac{(k-1)+(2k+2)\eta}{k-\frac{1}{16k^2}}\geq \frac{4}{5k}.
	\end{equation*}
	
	Now consider the matrix $B=A(G_{n,T} [E])$ and vector $y:= x_{|E}$ (the restriction of the vector $x$ to the set $E$). We see that for any vertex $v\in E$,
	$$
	By_v=\sum\limits_{\substack{u\sim v\\u\in E}}{x_u}\geq \frac{4}{5k}\rho x_v=\frac{4}{5k}\rho y_v.
	$$

	Hence, by Lemma \ref{lemma 3.1}, we have that $\rho(B)\geq \frac{4}{5k}\rho \geq \frac{4}{5}\sqrt{\frac{n-k}{k}}$. By Lemma \ref{lemma 3.4}, $\rho(B)\leq \sqrt{3k|E|}\leq \sqrt{3k\frac{n}{8k^2}}=\sqrt{\frac{3n}{8k}}$ when $|E|>k+2$. If $|E|\leq k+1$, since $G_{n,T}[E]$ is a bipartite graph, $\rho(B)=\rho(G_{n,T}[E])\leq \sqrt{\lfloor \frac{|E|}{2}\rfloor \lceil \frac{|E|}{2} \rceil }\leq\sqrt{\frac{|E|^2}{4}}<\sqrt{3k|E|}\leq \sqrt{\frac{3n}{8k}}$. By a simple calculation, we get a contradiction for large $n$. Thus, $E=\emptyset$. This implies that $V(G_{n,T})=L^{'}\cup R$. Note that $R$ is the common neighbourhood of $L^{'}$. Then $K_{k,n-k}\subset G_{n,T}$. Clearly, $G_{n,T}\cong K_{k,n-k}$. By the arbitrariness of $T$, we complete the proof of Theorem \ref{th4}. $\hfill\square$

	\section{The proof of Theorem \ref{th6}}\label{sec6}
	In this section, we give the proof of Theorem \ref{th6}. Before beginning our proof, we first give some notation not defined above. Denote by $F(G)$  the face set of $G,$ and let $\phi(G)=|F(G)|.$ For any face $f\in F(G)$, let $d_G(f)$ be the degree of the face $f$ in $G$. For short, $d_G(f)=d(f)$ if $G$ is clear from the context.
	
	\begin{lemma}\label{outerplanarEuler}
		If $G$ is a connected outerplanar graph of order $n$ and size $m$, then $m\leq 2n-3$. Moreover, if $G$ is a connected outerplanar bipartite graph, then $m\leq \frac{3n-4}{2}$.
	\end{lemma} 
	
	\begin{proof}
		Since $G$ is a connected outerplanar graph,   there exists a face in $G $ satisfying $d(f)=n.$ 
		Combining this with the fact that  $d(f)\ge 3$ for all $f\in F(G)$, we deduce that
		\begin{equation*}
			\sum_{f\in F(G)} d(f)\ge  3(\phi(G) -1)+n
		\end{equation*}
		Note that $\sum_{f\in F(G)} d(f)=2m$, and $n-m+\phi\ge 2$ by Euler's formula. Then by a simple calculation,  we yield that $m\le 2n-3$, as required.
		
		Furthermore, if $G$ is a  connected bipartite outerplanar graph, then $d(f)\ge 4$ for all $f\in F(G)$. Thus,
		$\sum_{f\in F(G)} d(f)\ge 4(\phi(G)-1)+n$. By a similar analysis as above, we obtain that
		$m\le \frac{3n-4}{2}$, as required.
	\end{proof}

	Now we shall give the proof of Theorem \ref{th6}.
	
	\noindent{\underline{\textbf{The proof of Theorem \ref{th6}. }}} Let $G$ be the outerplanar bipartite graph on $n$ vertices with maximum spectral radius $\rho(G)$.
	We first assume that $G$ is disconnected. Let $G_1,...,G_\eta$ be all components of $G$ for $\eta\ge 2$, and let
	$G'$ be a connected graph obtained from $G$ by adding $\eta-1$ edges.
	By the Rayleigh quotient and the Perron-Frobenius theorem,  $\rho(G')> \rho(G)$. Also note that $G'$ is also an outerplanar bipartite graph.  It is a contradiction because of the maximality of $\rho(G)$. For this reason,  in what follows, we always assume that $G$ is connected.
	
	Since the graph $K_{1,n-1}$ is an outerplanar bipartite graph, we obtain that
	\begin{equation}\label{1}
		\rho(G)\ge \rho(K_{1,n-1})= \sqrt{n-1}. 
	\end{equation}
	Let $x=(x_1,\cdots,x_n)$ be the eigenvector corresponding to the eigenvalue $\rho(G)$  such that $x_z=\max_{v\in V(G)} x_v =1.$ 
	In what follows, we shall divide our proof into three claims.
	
	\begin{claim}\label{claim1}
		For any vertex $u\in V(G)$,  $d(u)>x_un-7\sqrt{n}.$
	\end{claim}
	\renewcommand\proofname{\bf Proof}
	\begin{proof}
		Let $A_1=N_G(u)$ and $B_1=V(G)\setminus(A_1\cup \{u\}).$
		By eigenvalue-eigenvector equation, we have that
		\begin{equation}\label{2}
			\rho^2(G)x_u=\sum_{y\in N(u)}\sum_{v\in N(y)}x_v\le d(u)+\sum_{y\in N(u)}\sum_{v\in N(y)\cap A_1}x_v+\sum_{y\in N(u)}\sum_{v\in N(y)\cap B_1}x_v.
		\end{equation}
		Recall that $G$ is outerplanar bipartite graph. Then $e(A_1)=0$, and $G$ contain no $K_{2,3}$ as a subgraph which implies that each vertex in $B_1$ at most has two neighbors in $A_1$. Thus,
		\begin{equation}\label{3.0}
			\sum_{y\in N(u)}\sum_{v\in N(y)\cap A_1}x_v=0.
		\end{equation}
		By Lemma \ref{outerplanarEuler}(i.e. $e(G)\le\frac{3n-4}{2}$) 
		\begin{equation*}
			\begin{aligned}
			\sum\limits_{y\in N(u)}\sum\limits_{v\in N(y)\cap B_1}x_v
			\le 2\sum\limits_{v\in B_1} x_v
			\le \frac{2}{\rho(G)}\sum\limits_{v\in B_1}d(v)
			\le \frac{4e(G)}{\rho(G)}
			\le \frac{6n-8}{\rho(G)}
		\end{aligned}
	\end{equation*}
Combined this with (\ref{1}), 
\begin{equation}\label{3}
	\begin{aligned}
			\sum\limits_{y\in N(u)}\sum\limits_{v\in N(y)\cap B_1}x_v< \frac{6n-6}{\sqrt{n-1}}=6\sqrt{n-1}<6\sqrt{n},
			\end{aligned}
		\end{equation}
		According to (\ref{1}), (\ref{2}), (\ref{3.0}) and (\ref{3}), we obtain $\rho^2(G)x_u\le d(u)+6\sqrt{n}.$
		Again by (\ref{1}),  we yield that $d(u)>x_un-7\sqrt{n},$ as required.
	\end{proof}
	
	\begin{claim}\label{claim2}
		For any vertex $u\in V(G)\setminus\{z\}$ ,  $x_u <\frac{15}{\sqrt{n}}.$
	\end{claim}
	\renewcommand\proofname{\bf Proof}
	\begin{proof}
		Recall that $x_z=1$. Then by Claim \ref{claim1}, we have $d(z)> n-7\sqrt{n}.$ Thus, for any vertex $u\in V(G)\setminus\{z\}$, we have $d(u)<8\sqrt{n}.$ Since otherwise, $|N(u)\cap N(z)|>\sqrt{n}>3$, and so we can find a $K_{2,3}$ in $G$, a contradiction. Combining this with Claim \ref{claim1}, we have 
		$ 8\sqrt{n}>d(u)> x_un-7\sqrt{n},$ which follows $x_u <\frac{15}{\sqrt{n}},$ as required.
	\end{proof}
	
	\begin{claim}\label{claim3}
		$d(z)=n-1.$ 
	\end{claim}
	\renewcommand\proofname{\bf Proof}
	\begin{proof}
		Let $A_2=N(z)$ and $B_2=V(G)\setminus(A_2\cup \{z\}).$ Then $n=|A_2|+|B_2|+1.$
		We first prove that $\sum_{v\in B_2} x_v< \frac{557}{\sqrt{n}}.$ 
		In fact, combining (\ref{1}) with Claim \ref{claim2}, we have
		\begin{equation}\label{4}
			\sum_{v\in B_2}x_v \le \frac{ \sum_{v\in B_2}d(v)\cdot \frac{15}{\sqrt{n}}}{\rho(G)}\le \frac{15(e(A_2,B_2)+2e(B_2))}{\sqrt{n(n-1)}}.
		\end{equation}
		By Claim \ref{claim1}, we have $|A_2|=d(z)>n-7\sqrt{n}$, which follows that $|B_2|<7\sqrt{n}.$
		Since $G$  neither contain    $K_{2,3}$ nor $C_3$ as a subgraph,  we obtain that $e(A_2)=0$ and each vertex in $B_2$ at most has two neighbors in $A_2$, which follows $e(A_2,B_2)\le 2|B_2|<14 \sqrt{n}$. Combining this with Lemma \ref{outerplanarEuler}, we obtain $e(B_2)\le \frac{3|B_2|-4}{2}$, and so $2e(B_2)\le 3|B_2|-4<21\sqrt{n}.$
		Thus, (\ref{4}) becomes 
		\begin{equation}\label{5} 
			\sum_{v\in B_2}x_v < \frac{15\cdot 35}{\sqrt{n-1}}<\frac{557}{\sqrt{n}},
		\end{equation}
		for $n\ge 9,$ as required. 
		Next, we shall prove $B_2=\emptyset$ to complete our proof of the claim.
		Suppose to the contrary that $B_2\neq \emptyset$ and let vertex $y\in B_2$. 
		Recall that each vertex in $B_2$ at most has two neighbors in $A_2$. Combining this with Claim \ref{claim2} and (\ref{5}), we have that for $n>344569,$
		\begin{eqnarray*}
			\sum_{u\in N(y)}x_u &=& \sum_{u\in N(y)\cap B_2}x_u+\sum_{u\in N(y)\cap A_2}x_u\\
			&\le&  \sum_{u\in B_2}x_u +\frac{2\cdot 15}{\sqrt{n}} \\
			&<& \frac{(557+30)}{\sqrt{n}}<1.
		\end{eqnarray*}
		Let $G^*=G-\{yv:v\in N(y)\}+\{zy\}$. Then clearly, $G^*$ is also an outerplanar bipartite graph. However,
		$$\rho(G^*)-\rho(G)\ge \frac{x^\top(A(G^*)-A(G)x)}{x^\top x}=\frac{2x_y}{x^\top x}(1-\sum_{u\in N(y)}x_u)>0,$$
		which follows that $\rho(G^*)>\rho(G)$, contrary to the maximality of $\rho(G).$  Thus, $B_2=\emptyset$, as desired.
	\end{proof}
	Since $G$ is a bipartite graph, by Claim \ref{claim3}, we have $G\cong K_{1,n-1}.$ Thus, we complete the proof. $\hfill\square$

\end{document}